\documentclass[conference]{IEEEtran}
\IEEEoverridecommandlockouts
\usepackage{cite}
\usepackage{amsmath,amssymb,amsfonts}
\usepackage{algorithmic}
\usepackage{graphicx}
\usepackage{textcomp}
\usepackage{xcolor}
\usepackage{hyperref}
\usepackage{amsmath}
\usepackage{amsthm}
\usepackage{amssymb}
\usepackage{amsfonts}
\usepackage{xcolor}
\usepackage{mathtools}
\usepackage{verbatim}
\usepackage{multicol}
\usepackage{multirow}
\usepackage{mathabx}
\usepackage{graphicx}
\usepackage{booktabs}
\usepackage{float}
\usepackage{listings}
\usepackage{array}
\usepackage{etoolbox}
\usepackage{setspace}
\usepackage{hyperref}
\hypersetup{
    colorlinks=black,
    linkcolor=black,
    filecolor=black,      
    urlcolor=black,
    citecolor=black
}
\usepackage{subcaption}
\usepackage{xr}
\usepackage{array}
\makeatletter
\def\th@plain{%
  \thm@notefont{}
  \itshape 
}
\def\th@definition{%
  \thm@notefont{}
  \normalfont 
}

\makeatother
\theoremstyle{plain}
   \newtheorem{theorem}{Theorem}[section]
   \newtheorem*{theorem*}{Theorem}
   \newtheorem{proposition}[theorem]{Proposition}

   \newtheorem*{lemma*}{Lemma}
   \newtheorem*{corollary*}{Corollary}
\theoremstyle{definition}
   \newtheorem{definition}{Definition}[section]

\definecolor{green}{RGB}{34, 139, 34}
 \numberwithin{equation}{section}

\theoremstyle{remark}

\newtheorem{remark}{Remark}[section]

\newcommand{\RR}{\mathbb{R}}

\def\BibTeX{{\rm B\kern-.05em{\sc i\kern-.025em b}\kern-.08em
    T\kern-.1667em\lower.7ex\hbox{E}\kern-.125emX}}
\begin{document}

\title{Neumann eigenmaps for landmark embedding\\
}

\author{\IEEEauthorblockN{Shashank Sule}
\IEEEauthorblockA{\textit{Department of Mathematics} \\
\textit{University of Maryland, College Park}\\
College Park, United States of America \\
ssule25@umd.edu}
\and
\IEEEauthorblockN{Wojciech Czaja}
\IEEEauthorblockA{\textit{Department of Mathematics} \\
\textit{University of Maryland, College Park}\\
College Park, United States of America \\
wojtek@math.umd.edu}

}

\maketitle

\begin{abstract}
We present Neumann eigenmaps (NeuMaps), a novel approach for enhancing the standard diffusion map embedding using \emph{landmarks}—distinguished samples within the dataset. By interpreting these landmarks as a subgraph of the larger data graph, NeuMaps are obtained via the eigendecomposition of a renormalized Neumann Laplacian. We show that NeuMaps offer two key advantages: (1) they provide a computationally efficient embedding that accurately recovers the diffusion distance associated with the reflecting random walk on the subgraph, and (2) they naturally incorporate the Nystr\"om extension within the diffusion map framework through the discrete Neumann boundary condition. Through examples in digit classification and molecular dynamics, we demonstrate that NeuMaps not only improve upon existing landmark-based embedding methods but also enhance the stability of diffusion map embeddings to the removal of highly significant points.
\end{abstract}

\begin{IEEEkeywords}
Manifold learning, Diffusion maps, Spectral graph theory.
\end{IEEEkeywords}

\section{Introduction}
Manifold learning algorithms, such as principal component analysis, locally linear embeddings, diffusion maps, and their variants, are essential tools for unsupervised learning in sophisticated real-world datasets. However, the low-dimensional embeddings produced by these algorithms can become computationally intractable due to an eigendecomposition or alignment step which scales poorly with sample complexity. To address this challenge, a variety of \emph{landmarking} methods have been developed to reduce computational load and to incorporate distinguished data points into the unsupervised learning process \cite{landmark1, landmark2, landmark3, landmark4, landmark5, landmark6}.  

Typically, in landmarking methods, a subset of landmarks $ V_S $ is selected from the dataset $V_G := \{x_i\}_{i=1}^{n} \subseteq \mathbb{R}^m $, where $|V_S| \ll |V_G| = n$. These landmarks can be chosen either randomly or based on their significance within the dataset. The dimensionality reduction technique is then efficiently applied to $ V_S $, resulting in an embedding $ \psi: V_S \to \mathbb{R}^d $, where $ d < m $. Subsequently, $\psi $ is extended to the full dataset $ V_G $ using a computationally inexpensive out-of-sample extension algorithm. This raises the natural question of \emph{how to integrate the rest of the data $\delta S:= V_G \setminus V_S$ into the dimensionality reduction process while retaining the computational speedup of focusing on $V_S$}. This question is relevant even when computational complexity is not an issue: in particular, embeddings of $V_S$ that account for $\delta S$ may provide advantages over embeddings that do not. In this context, we may even \emph{switch} the roles of the landmarks: for instance, in \cite{landmarkdmap2} it was suggested that landmarks be $\delta S$ instead and inform the embedding of the rest of the set $V_S$. 

In this paper, we address this question specifically for the diffusion map (Dmap) algorithm by introducing Neumann maps (NeuMaps). In particular, by envisioning the overall data $V_G$ as the vertices of a graph, we propose that the subset $V_S \subseteq V_G$ be embedded via the normalized Neumann eigenvectors of the subgraph $S$ induced by $V_S$. The Neumann boundary condition generates the reflected random walk on subgraphs where the random walker may diffuse between landmarks either directly or via a reflection off the boundary vertices in $\delta S = V_G \setminus V_S$. The diffusion distance so defined by this random walk is thus additionally robust to perturbations due to data subsampling and emphasizes cluster structure due to the added within-cluster diffusion probability via reflection. 

\subsection{Related work} The Neumann spectrum of subgraphs and the Neumann Laplacian were established in the seminal works by Chung, Graham, and Yau \cite{sgt1, sgt2, sgt3} and later chronicled in the monograph by Chung \cite{SGTbook}. Here we study the properties of these eigenvectors as landmark embeddings. Integrating landmarks within the diffusion maps framework--and in particular using the removed data to inform landmark embeddings--is an active topic of research in manifold learning \cite{landmark0, landmarkdmap1, landmarkdmap2, landmarkdmap3, landmarkdmap4, landmarkdmap5}. 

\subsection{Contributions and organization} In this paper, we demonstrate the utility of the normalized Neumann eigenvectors of subgraphs as a landmark-based manifold learning technique. To this end, in Section 2 we briefly review the theory of Neumann eigenvectors of subgraphs from \cite{SGTbook}. In Section 3 we define Neumann eigenmaps as the eigenvectors of the reflecting random walk transition kernel; this definition allows us to prove an isometric embedding result for recovering a probability distance metric on the data. In Section 4 we provide numerical examples that benchmark NeuMaps against \emph{Roseland} \cite{landmark2} and Dmaps.

\section{Neumann spectra of subgraphs}\label{sec: Neumann and Dirichlet Spectra}

We will consider the standard setup in diffusion maps where the data $V_G$ is used to construct a weighted graph where the edge weights are described by an affinity kernel. We therefore switch to using terminology on graphs. Moreover, note that if $A \in \mathbb{R}^{n\times m}$ and $X \subseteq [n], Y\subseteq [m]$ then $A[X,Y]$ denotes the submatrix obtained by selecting rows w.r.t in $X$ and columns w.r.t indices in $Y$. If $X = \{x\}, Y = \{y\}$ then $A[X,Y] = A[x,y]$, the $xy$\,th entry in $A$. 

\begin{definition}\label{def: graph}
Let $G = (V_G, E_G, W_G)$ be a finite weighted graph with vertices $V, |V| = n$, an adjacency matrix $W_{G} \in \mathbb{R}^{n \times n}$. Given an enumeration of $V_G$ we identify it with $[n]$. Then the edges are defined as $E_{G} = \{\{{x,y}\} \mid x, y \in V_G, w(x,y): = W[x,y] > 0\}$. 
\end{definition}
\begin{definition}\label{def: degree and laplacian}
    Let $\mathbf{1}_{|V_G|}$ be the vector of 1's in $|V_G|$ dimensions. The degree matrix of $G$ is given by $D_G = \textsf{diag}(W_G \mathbf{1})$ and the graph Laplacian is given by $L_G = D_G - W_G$. For any $x \in V_G$, we denote $d(x) := D_G[x,x]$. 
\end{definition}
\begin{definition}\label{def: subgraph info}
    Let $V_S \subseteq V_G$. Identifying $V_S$ with the corresponding subset of $[n]$, the subgraph induced by $V_S$ is given by $S = (V_S, E_{S}, W_S)$.
\end{definition}
\begin{definition}\label{def: degree info}
Let $V_S \subseteq V_G$ and $S = (V_S, E_S, W_S)$. The \emph{graph degree matrix} of $S$ is given by $T_{S} = D_G[V_S,V_S]$. Note that for any $x \in V_S$, we have $T_S[x,x] = d(x)$. 
\end{definition}

\begin{definition}\label{def: boundary info}
Let $S$ be a subgraph of $G$. We define the boundary edges $\partial S$ and boundary vertices $\delta S$ as: 
\begin{align}
    \partial S &= \{\{x,y\} \mid x \in V_S, y \in  V_G \setminus V_S, w(x,y) > 0\}, \\
    \delta S &= \{y \in V_G \setminus V_S \mid \exists \: e \in \partial S \text{ s.t. } y \in e\}. 
\end{align}
Moreover let $S^* = \partial S \cup E_S$. 
\end{definition}
\begin{definition}\label{ch2_def: Neumann eigenvalues}
Let $S \subseteq G$. Then the \textit{Neumann eigenfunction} $f^{N}_{1}: S \cup \delta S \to \mathbb{R}$ and $\textit{Neumann eigenvalue}$ $\lambda_{1}^{N}$ of $S$ are defined as follows: 
\begin{align}
    \lambda^{N}_1 &= \underset{f \mid_S\perp T_{S}1}{\text{min}} \frac{\sum_{\{x,y\} \in S^*}w(x,y)(f(x) - f(y))^2}{\sum_{x \in V_S}(f(x))^2 d(x)}\label{eq: Neumann eigenvalue}.
\end{align}
The minimizer in \eqref{eq: Neumann eigenvalue} will be termed $f^{N}_1$. In general, we may sequentially generate the $i$th Neumann eigenpair by constraining the minimization problem \eqref{eq: Neumann eigenvalue} the orthogonal subspace of $\textsf{span}\{f_1, \ldots, f_{i-1} \}.$
\end{definition}



The following result proves that the Neumann eigenfunction satisfies a vanishing discrete normal derivative on the boundary $\delta S$:

\begin{theorem}[Lemma 8.1 in \cite{SGTbook}]\label{thm: Neummann eigenfunction main theorem}
Let $V_S \subseteq V_G$ and $L_{V^*}$ the Laplacian of the graph generated by $V^*=V_S \cup \delta S$. The Neumann eigenfunction $f = f_{1}^{N}$ satisfies the following properties: 
\begin{enumerate}
    \item Fix $x \in V_S$. Then 
    \begin{equation}
        L_{V^*}f(x) = \lambda_{N}^{1}d(x)f(x).\label{eq: Eigenvalue equation Neumann}
    \end{equation}
    \item Fix $x \in \delta S$. Then 
    \begin{equation}
        \sum_{\{x,y\} \in \partial S}w(x,y)(f(x) - f(y)) = 0. \label{eq: Neumann condition}
    \end{equation}
\end{enumerate}
\end{theorem}

\begin{remark}
The condition \eqref{eq: Eigenvalue equation Neumann} shows that $f^{N}_1$ satisfies an eigenvalue equation at $x \in V_S$ for the \emph{ambient graph} $G$, while not necessarily being either a Laplacian eigenvector of the graph $S$ or $G$. Moreover, condition \eqref{eq: Neumann condition} can be viewed as a vanishing discrete normal derivative condition. 
\end{remark}

\section{Neumann Maps}\label{sec: Neumann maps}
The quotient \eqref{eq: Neumann eigenvalue} modifies the usual normalized Rayleigh quotient of the graph Laplacian on $S$ by penalizing large variations between $V_S$ and $\delta S$, thus motivating the use of $f^{N}_i$ as feature maps specifically adapted for smooth extensions to $\delta S$. We now show that the minimizers of \eqref{eq: Neumann eigenvalue} can be computed via the eigenvectors of the \emph{Neumann Laplacian} $L^{N}_{S}$. First, we define additional terms below. 
\begin{definition}
    The \emph{boundary operator} is defined as $B_S = W_{G}[V_{G}\setminus V_S, V_S]$ and the boundary degree matrix is defined by $T^{\delta}_{S} = \textsf{diag}(B_{S}\mathbf{1}_{|V_S|}) \in \mathbb{R}^{|\delta S| \times |\delta S|}$. Here $\mathbf{1}_{|V_S|} \in \mathbb{R}^{|V_S|}$. 
\end{definition}
The boundary operator can be used to define a \emph{Neumann Laplacian}:
\begin{definition}
     Let $S= (V_S, E_S, W_S)$ be a subgraph of $G$ and $L^{D}_{S} := L_{G}[V_S,V_S]$. Moreover, let $B_S$ be the boundary matrix. Then the Neumann Laplacian, $L^{N}_{S}$, is defined as follows: 
     \begin{align}
          L^{N}_{S} := L^{D}_{S} - B^{\top}(T^{\delta}_{S})^{-1}B.
     \end{align}
\end{definition}
The following proposition demonstrates how the Neumann Laplacian reformulates conditions \eqref{eq: Neumann condition} and \eqref{eq: Eigenvalue equation Neumann} into one eigenproblem: 
\begin{proposition}
\label{prop: Neumann laplacian}
Let $v: V_S \cup \delta S \to \mathbb{R}$ satisfy the Neumann condition \eqref{eq: Neumann condition} on $\delta S$ and $u = v\mid_{V_S}$ be its restriction to $S$. Then $v$ satisfies the Laplacian eigenvalue condition \eqref{eq: Eigenvalue equation Neumann} if and only if $u$ satisfies: 
    \begin{align}
        \lambda_{1}^{N}T_{S}u = L^{N}_{S}u. \label{eq: Neumann Laplacian eigenproblem}
    \end{align}
\end{proposition}
One may thus obtain the Neumann eigenvectors by computing an eigenvector of $L^{N}_{S}$ and then extending it to $\delta S$ via the Neumann condition \eqref{eq: Neumann condition}. Then from the converse of Proposition \eqref{prop: Neumann laplacian} the resulting vector on $V_S \cup \delta S$ is a Neumann eigenfunction. More importantly, the renormalized Neumann Laplacian $\mathcal{N} := T_{S}^{-1/2}L^{N}_{S}T_{S}^{-1/2}$ is similar to a random walk matrix:

\begin{proposition}\label{prop: Neumann random walk matrix is a random walk matrix}
Let $R = T_{S}^{-1/2}(I - \mathcal{N})T_{S}^{1/2}$. Then $R\textbf{1} = \textbf{1}$ and the entries of $R$ are all non-negative. 
\end{proposition}

We can now define a Neumann map by diagonalizing $R$ via the spectral decomposition of $\mathcal{N}$:  $$R = T_{S}^{-1/2}(I - \mathcal{N})T_{S}^{1/2} = T_{S}^{-1/2}U\Sigma U^{\top}T_{S}^{1/2} = V\Sigma Y^{\top}$$
where $V = T_{S}^{-1/2}U$ and $Y = T_{S}^{1/2}U$. Since $Y^{\top}V = I$, for any $t \geq 0$ we have $$R^{t} = \sum_{i=1}^{n}\sigma_{i}^{t}v_{i}y^{\top}_{i}$$ so in particular expanding the column space of $R^t$ we get
\begin{align}
   R^{t}[:,x_i] = \sum_{j=1}^{n}\sigma_{j}^{t}v_j y_{j}^{\top}(x_i). \label{eq: column space decomposition} 
\end{align}

\begin{definition}\label{ch3_def: Neumann diffusions}
Consider the column space decomposition of $R^{t}[:,i]$ in \eqref{eq: column space decomposition}. The $d$-dimensional Neumann map of vertex $i \in S$ is the following point in $\RR^{d}$:
\begin{align*} g^{d}_{t}(x_i) = \begin{bmatrix}\sigma_{2}^{t}y_{2}^{\top}(x_i) \\ \vdots \\ \sigma_{d+1}^{t}y_{d+1}^{\top}(x_i)\end{bmatrix}.\end{align*}
\end{definition}

\begin{theorem}\label{thm: Main theorem of Neumann diffusions} 
Let $X_t$ be the random walk defined by $R$ and $p^{t}_{ji} := P(X(t) = i \mid X(0) = j)$ be the probability of walking to vertex $i$ after starting at vertex $j$ after $t$ steps. Then
\begin{align}
&\sum_{j=1}^{|S|}(p^{t}_{ji} - p^{t}_{jk})d(x_j)^{-1} = ||g^{|S|}_{t}(x_{i}) - g^{|S|}_{t}(x_{k})||^2 .
\end{align}

\end{theorem}

\begin{remark}\label{ch3_rem: reflecting walks are neumann walks} The matrix $R$ is the transition matrix of the reflecting random walk mentioned in \cite{SGTbook} wherein if $u,v \in S$ then one may walk from $u$ to $v$ either directly or via a reflection off some mutually adjacent $x \in \delta S$. This random walk may be contrasted to the one in Roseland which considers random walks between vertices in $S$ \emph{only} via walks through $\delta S$. Therefore, the choice to include within-$S$ diffusion defines the conceptual difference between NeuMaps and Roseland. Moreover, note that $L^{N}_{S}$ is a perturbation of $L^{D}_{S}$, the Dirichlet graph Laplacian (see \eqref{eq: Neumann Laplacian eigenproblem}). Therefore, in theory, alternate Neumann Laplacians could be obtained by replacing $L^{D}_{S}$ with any user-chosen operator such as the Schr\"odinger or transport operator on the subgraph $S$ \cite{landmark0, landmarkdmap4, landmarkdmap5}. 
\end{remark}


\begin{remark} The Neumann condition gives a natural way to extend functions from $V_S$ to $V_S \cup \delta S$. In particular, if $\mathcal{N}g = \lambda g$ then defining $f = T^{1/2}g$ we have $L^{N}_{S}g = \lambda T_{S}g$. To make $f$ a Neumann eigenvector, we set its discrete normal derivative to zero, which corresponds to setting $f(x) = d(x)^{-1}\sum_{y \in V_S}w(x,y)f(y)$ for $x \in \delta S$. By rescaling, we get the following formula for extending $g$: 
    \begin{align}
        g(x) = \sum_{y \in V_S}\frac{w(x,y)}{\sqrt{d(x)d(y)}}g(y). \label{eq: extension formula}
    \end{align}
The extension of $g$ via \eqref{eq: extension formula} is therefore a multiple of its Nystr\"om extension \cite{nystrom}. Consequently, the Neumann condition naturally accounts for Nystr\"om extension, giving an additional justification for using NeuMap. 
\end{remark}


\section{Numerical experiments}

In the examples below, we adopt the paradigm where the landmarks are $\delta S$ and the rest of the set $V_S$ is embedded using NeuMap. We remark that in both examples below, $|\delta S| \leq |V_G \setminus \delta S|$. Thus, the computational acceleration of using NeuMap instead of Dmap will be modest. However, here our goal is not necessarily to emphasize the computational efficiency but the more favourable geometric and structural properties of the NeuMap embedding. 
\subsubsection{Digit classification and Roseland}
\begin{figure*}[t]
    \centering
    \includegraphics[width=0.7\textwidth]{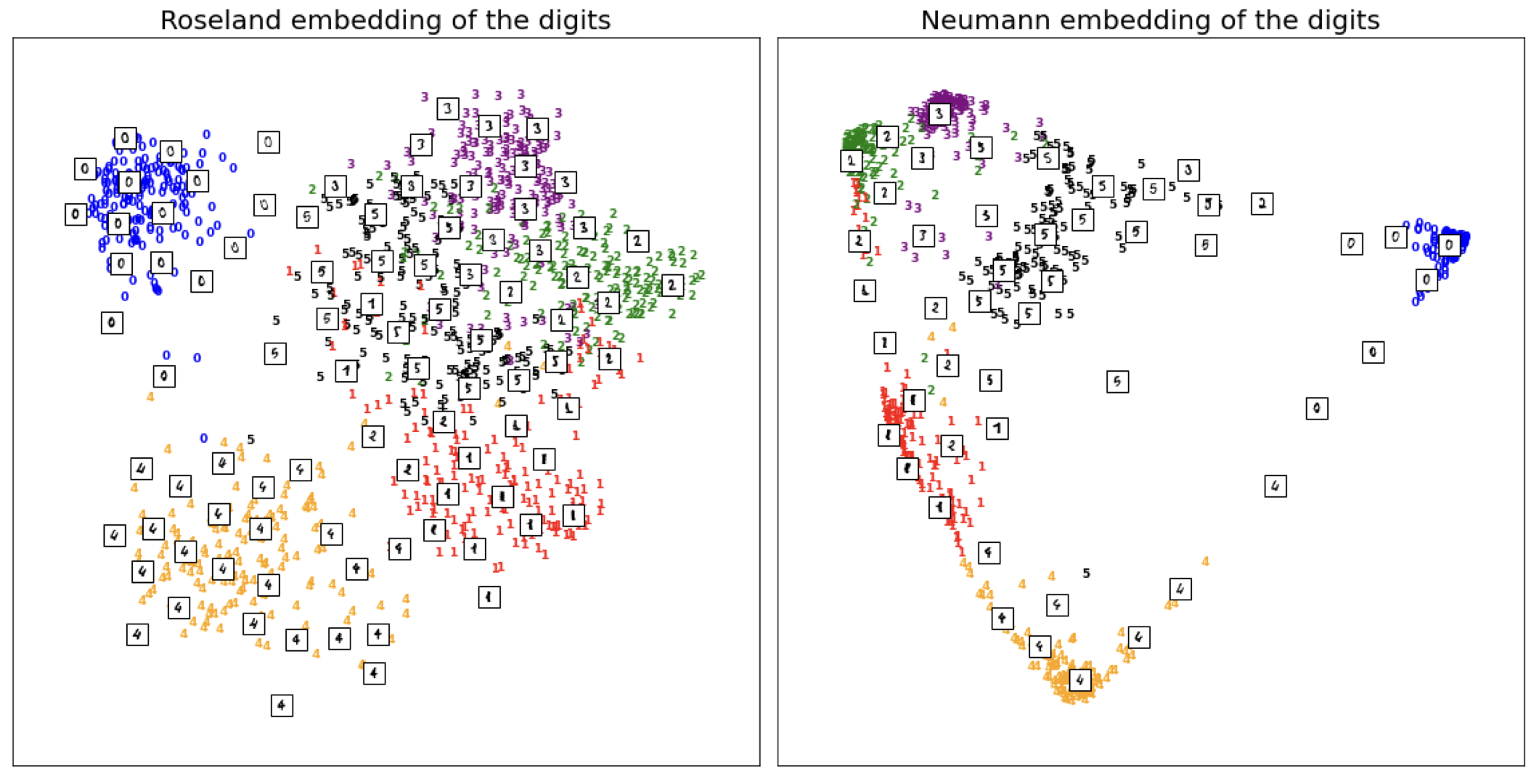}
    \caption{Left: Embedding UCI digits with 25\% random landmarking via Roseland. Right: Embedding UCI digits with Neumann maps.}
    \label{fig: digits}
\end{figure*}

We consider the UCI digits dataset of 1083 $16 \times 16$ images of handwritten digits from $0-6$. We randomly select 25\% of the data as landmarks. As a comparison, we use the Roseland algorithm. We tune the bandwidth $\epsilon$ according to the max-min criterion in \cite{datafold}. To keep the comparisons consistent, we use the same $\epsilon$ for computing Neumann maps. After projecting in 2D we assign clusters using $k$-means and measure the cluster assignments using the normalized mutual information (NMI) and  unsupervised clustering accuracy (ACC) against the true labels. For Roseland, the NMI and ACC were 0.71 and  84\% respectively while for Neumann maps the NMI and ACC were 0.85 and 93\% respectively. Thus, Neumann maps outperforms Roseland on this task. We further visualize our results in Figure \ref{fig: digits} where we observe that Roseland by and large successfully separates the data in 6 clusters. However, in Neumann maps the clusters are significantly more concentrated, likely contributing to the higher accuracy of $k$-means on this embedding. An interesting artefact of the Neumann embedding seems to be the presence of few but significant outliers which are rather far apart from the cluster mean. 

\subsubsection{Learning collective variables in molecular dynamics}

In the analysis of data from molecular dynamics (MD) simulations, a highly important task is the recovery of optimal \emph{collective variables}, i.e low dimensional features of the high dimensional atomic coordinates \cite{cvsimportant}. Criteria for optimality of such CVs specifically geared for MD include accurate representations of the dynamical properties of the system \cite{cvsdynamical1, cvsdynamical2, cvsdynamical3}, reproducibility of transition rates and reaction channels in low dimensions \cite{cvsrr1, cvsrr2}, and ability to drive enhanced sampling simulations. 

\begin{figure}
    \centering
    \includegraphics[width=.9\linewidth]{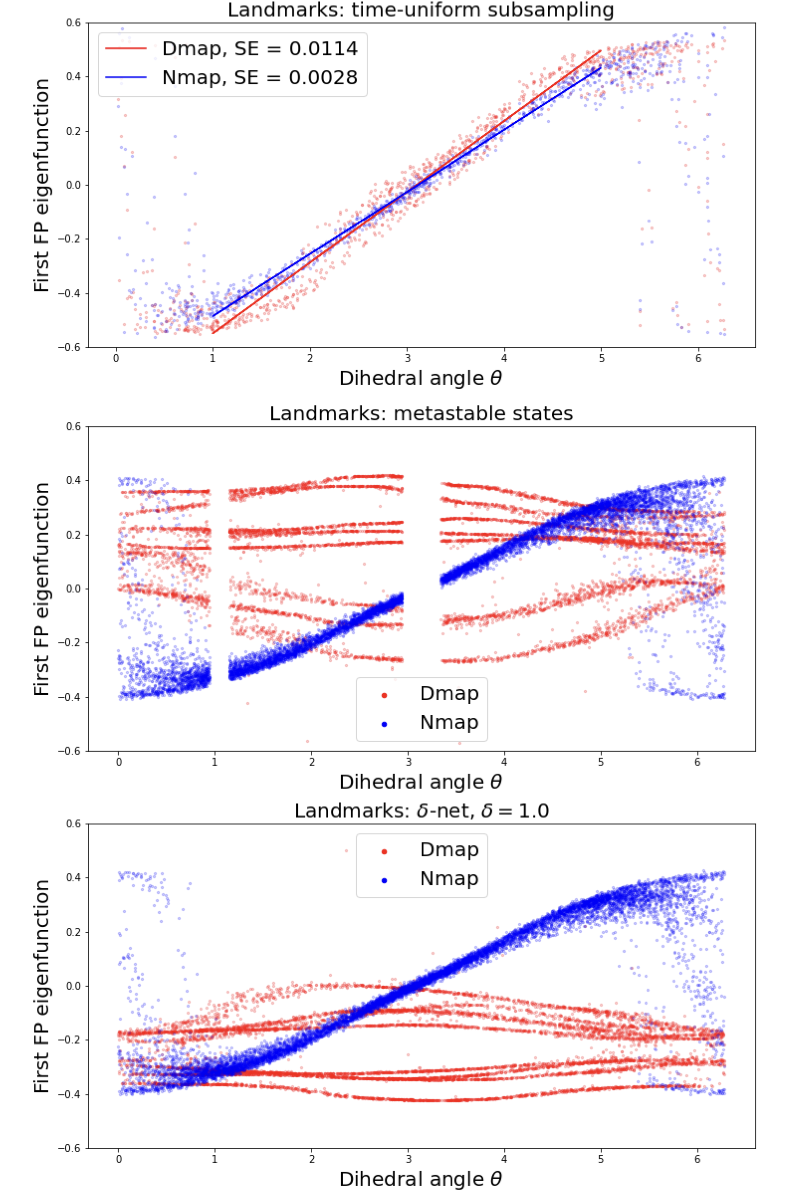}
    \caption{Top to bottom: FP eigenfunction embeddings of data after removal of landmarks provided by subsampling uniformly in time, choosing metastable states, and uniformly subsampling in space via $\delta$-nets.}
    \label{fig: butane}
\end{figure}

While several algorithms for learning CVs have been proposed \cite{cv1, cv2, cv3, cv4, cv5}, the Fokker-Planck eigenfunctions--computable via Dmap--stand out due to their simplicity and their optimal embedding properties in a probability distance metric \cite{fp1, fp2}. We now show that given certain marked points, the Neumann Fokker-Planck eigenfunctions can also be used to learn CVs and visualize MD data. In particular, we demonstrate that they may learn the underlying collective variable \emph{more accurately} than Fokker-Planck (FP) eigenfunctions using diffusion maps. We empirically illustrate this in the case of the butane ($C_4 H_{10}$) molecule, widely used as a toy model for configurational changes in small molecules. In particular, it is well-known that the dynamics of the butane molecule, while residing in a 42-dimensional space, are accurately coarse-grained by the one-dimensional \emph{dihedral angle} $\theta$ in its carbon backbone. The dihedral angle characterizes the two metastable states of the butane molecule: the anti configuration given by $\theta \approx \pi$ and the gauche configurations given by $\theta \approx \pi/3, 5\pi/3$. The atomic coordinates $X_t \in \mathbb{R}^{14 \times 3}$ largely inhabit the metastable states and transition rarely between them. Additional enhanced sampling algorithms such as metadynamics may be used to generate more samples in the transition regions. To this end, we simulate a metadynamics trajectory for $X_t$ in OpenMM using the Langevin integrator at 300K and collect $10^4$ points given by $\{X_{i\Delta t}\}_{i=1}^{10^4}$ where the integrator timestep $\Delta t = 0.04 $ps. In order to apply Neumann maps to this data, we consider the following ways to obtain marked and unmarked data: (1) Subsampling uniformly in time (i.e every 10th point is labeled as marked; the rest are unmarked), (2) Marking metastable states i.e points corresponding to $|\theta - \pi| \leq 0.2$ and $|\theta - \pi/3| < 0.1, |\theta - 5\pi/3| < 0.1$ are labeled as marked, and (3) Subsampling uniformly in ambient space, where we select a $\delta$-net such that all points are within distance $\delta$ of this net. In all three cases, approximately 10\% of the points are marked. We detail our results in Figure \ref{fig: butane}. In the top panel, we conisder the subsampling-uniformly-in-time scheme. The first FP eigenfunctions $\psi_1$ computed using Dmap and NeuMap both correlate significantly with $\theta$. To measure this correlation, we model $\theta$ as a linear function of $\psi_1$ and observe a standard error (SE) of $2.8 \times 10^{-4}$ for $\psi_1$ computed via NeuMap against a much higher SE of $1.1 \times 10^{-2}$ for $\psi_1$ computed via Dmap, thus showing that NeuMap recovers the dihedral angle more accurately. By removing metastable states (middle) and removing $\delta$-nets (bottom)--we find that the performance of the Dmap degrades significantly, thus revealing that these points were instrumental in producing the original embedding that recovered the dihedral angle. However, when computing with NeuMap--presumably due to the robustness provided by reflection against the removed points--the embedding stays relatively unperturbed. Therefore, by incorporating deleted points yet retaining the computational efficiency of using fewer samples, NeuMap makes Dmap stable to erasures. 

\section{Conclusion}
We have introduced NeuMap which enhances the Dmap embedding with landmarks and the Neumann boundary condition. The embedding so obtained recovers the diffusion distance of the reflecting random walk, accelerates the Dmap process via the eigendecomposition of a smaller matrix, and exhibits favourable performance over Roseland and Dmap in proof-of-concept examples.

\section*{Appendix}
Here we state the proofs of the results stated in the main paper. 

\subsection{Proof of theorem \eqref{thm: Neummann eigenfunction main theorem}}
We prove parts (1) and (2) variationally. For part (1),  fix $x_0 \in V_S$, define $vol(S) := \sum_{x \in V_S}d(x)$ and let 
    \begin{align*}f_{\epsilon}(x) = \begin{cases}f(x_0) + \frac{\epsilon}{ d_{x_0}} \text{ if } x=x_0, \\ f(x) - \frac{\epsilon}{ vol(S) - d_{x_0}}. \end{cases}\end{align*}
    First observe that $f_{\epsilon}\mid_{S} \perp T1_S$ so the minimization problem is well-defined on $f_{\epsilon}$. Now we compute the quotient for $f_{\epsilon}$:
    \begin{gather*}
        R(\epsilon) = \frac{\sum\limits_{(x,y) \in S^{*}}w(x,y)(f_{\epsilon}(x) - f_{\epsilon}(y))^2}{\sum\limits_{x \in V_S}(f_{\epsilon}(x))^2 d_x}\\
        = \frac{\sum\limits_{(x,y) \in S^{*}, x \neq x_0}w(x,y)(f(x) - f(y))^2}{\sum\limits_{x \neq x_0}(f(x) - \frac{\epsilon}{\text{vol}(S) - d_{x_0}})^2d_{x_0} + (f(x_0) + \frac{\epsilon}{d_{x_0}})^2d_{x_0}}\\
        + \frac{\sum\limits_{(x_0,y) \in S^{*}}w(x_0,y)(f(x) +  \frac{\epsilon}{d_{x_0}} - f(y) + \frac{\epsilon}{\text{vol}(S) - d_{x_0}}))^2}{\sum\limits_{x \neq x_0}(f(x) - \frac{\epsilon}{\text{vol}(S) - d_{x_0}})^2d_{x_0} + (f(x_0) + \frac{\epsilon}{d_{x_0}})^2d_{x_0}} \\
         = \frac{\sum\limits_{(x,y) \in S^{*}}w(x,y)(f(x) - f(y))^2}{\sum\limits_{x \in V_S}(f(x))^2 d_{x} + \frac{2\epsilon f(x_0)d_{x_0}\text{vol}(S)}{d_{x_0}(vol(S) - d_{x_0})} + O(\epsilon^2)}. \\
        + \frac{ \frac{2\epsilon\text{vol}(S)}{d_{x_0}(\text{vol}(S) - d_{x_0})}\sum\limits_{(x_0,y) \in S^*}w(x_0,y)(f(x_0) - f(y)) + O(\epsilon^2)}{\sum\limits_{x \in V_S}(f(x))^2 d_{x} + \frac{2\epsilon f(x_0)d_{x_0}\text{vol}(S)}{d_{x_0}(vol(S) - d_{x_0})} + O(\epsilon^2)}.
    \end{gather*}
The second equality follows after simplifying the algebra and noting that $\sum\limits_{x \in V_S}f(x)d_{x} = 0$. We know that when $\epsilon = 0$, $f_{\epsilon} = f$, which also minimizes $R(\epsilon)$. Thus, $R'(0) = 0$ so computing the derivative via the quotient rule and setting the numerator at $\epsilon = 0$ to zero, we get that 
\begin{align}
    &\Big(\frac{2\text{vol}(S)}{d_{x_0}(\text{vol}(S) - d_{x_0})}\cdot \\
    &\sum\limits_{(x_0,y) \in S^*}w(x_0,y)(f(x_0) - f(y))\Big)\sum\limits_{x \in V_S}(f(x))^2 d_{x}\\
    &- \Big(\frac{2f(x_0)d_{x_0}\text{vol}(S)}{d_{x_0}(\text{vol}(S) - d_{x_0})}\Big)\sum\limits_{(x,y) \in S^{*}}w(x,y)(f(x) - f(y))^2 = 0.
\end{align}
Rearranging the equation, dividing through by $\sum\limits_{x \in V_S}(f(x))^2 d_{x}$ and noting that 
$$\frac{\sum\limits_{(x,y) \in S^{*}}(f(x) - f(y))^2}{\sum\limits_{x \in V_S}(f(x))^2 d_{x}} = \lambda_{1}^{N}$$

gives us part (1). For part (2), we adopt a similar strategy but the variation is simpler. Fix $x_0 \in \delta S$ and set 
\begin{align*} f_{\epsilon}(x) = \begin{cases} f(x) + \epsilon \text{ if } x = x_0, \\ f(x) \text{ otherwise. } \end{cases}\end{align*}
Now compute the Neumann quotient and observe that we can separate the sum in the numerator over edges that connect to $x_0$ and those that don't. By definition, the edges that connect with $x_0$ are contained in $\partial S$ so 
\begin{align*}
R(\epsilon) &= \frac{\sum\limits_{(x,y) \in S^{*}}w(x,y)(f_{\epsilon}(x) - f_{\epsilon}(y))^2}{\sum\limits_{x \in V_S}(f_{\epsilon}(x))^2 d_x} \\
&= \frac{\sum\limits_{(x,y) \in S^{*}}w(x,y)(f(x) - f(y))^2}{\sum\limits_{x \in V_S}(f(x))^2 d_x} \\
&+ \frac{2\epsilon\sum\limits_{(x_0,y) \in \partial S}w(x_0,y)(f(x_0) - f(y)) + O(\epsilon^2)}{\sum\limits_{x \in V_S}(f(x))^2 d_x}.
\end{align*}
Once again, taking the derivative with respect to $\epsilon$ and setting it $0$ at $\epsilon = 0$ yields (2). 
\subsection{Proof of proposition \ref{prop: Neumann laplacian}}
Note that $v$ satisfies the Neumann condition if and only if for every $y \in \delta S$, 
\begin{align}
    f(y) = \frac{1}{\partial d(y)}\sum_{z \in V_S}w(y,z)v(z). \label{eq: neumann condition rephrased as boundary average}
\end{align} Here $\partial d(y) = T^{\delta}_{S}[y,y]$. To prove the equivalence stated, we compute the action of $L_{V^*}$ on $v$: 
\begin{align}
    &L_{V^*}v(x) = d(x)v(x) - \sum_{y \in V^*}w(x,y)v(y) \\
    &= d(x)v(x) - \sum_{y \in V}w(x,y)v(y) - \sum_{y \in \delta S}w(x,y)v(y) \\
    &= d(x)v(x) - \sum_{y \in V}w(x,y)v(y) - \sum_{y \in \delta S}\frac{w(x,y)}{\partial d(y)}\sum_{z \in V_S}w(y,z)v(z) \\
    &= d(x)u(x) - \sum_{y \in V}w(x,y)u(y) - \sum_{y \in \delta S}\frac{w(x,y)}{\partial d(y)}\sum_{z \in V_S}w(y,z)u(z) \\
    &= L_{S}^{N}u(x). 
\end{align}
Here the third equality follows by plugging in \eqref{eq: neumann condition rephrased as boundary average} and the fourth equality follows by noting that $v\mid_{V_S} = u$. Clearly, $L_{V^*}v(x) = \lambda d(x)v(x)$ if and only if $L_{S}^{N}u(x) = \lambda d(x)v(x)$. This proves our assertion. 

\subsection{Proof of proposition \ref{prop: Neumann random walk matrix is a random walk matrix}}

This follows because $\mathcal{N}$ admits 
$T^{1/2}\textbf{1}$ as a zero-eigenvector:

\begin{align*}R\textbf{1} = T_{S}^{-1/2}(I - \mathcal{N})T_{S}^{1/2}\textbf{1} = I\textbf{1} - T_{S}^{-1/2}\mathcal{N}T_{S}^{1/2}\textbf{1} = \textbf{1}.\end{align*}
To see the non-negativity of the entries, we expand $N$ in terms of the Dirichlet and Boundary operators:
\begin{align*}
    R &= T_{S}^{-1/2}(I - \mathcal{N})T_{S}^{1/2}\\
    &= I - T_{S}^{-1/2}\mathcal{N}T_{S}^{1/2}\\
    &= I - T_{S}^{-1}(L_G[V_S, V_S] - B^{\top}(T^{\delta}_{S})^{-1}BT_{S}^{-1/2})\\
    &= I - T_{S}^{-1}L_G[V_S, V_S] + T_{S}^{-1}B^{\top}(T^{\delta}_{S})^{-1}B \\
    &= I - T_{S}^{-1}(T_{S} - W_G[V_S, V_S]) + T_{S}^{-1}B^{\top}(T^{\delta}_{S})^{-1}B \\
    &= T_{S}^{-1}W_G[V_S, V_S] + T_{S}^{-1}B^{\top}(T^{\delta}_{S})^{-1}B
\end{align*}
Clearly the above matrix has nonnegative entries. Thus $R$ is row-stochastic. 
\subsection{Proof of \ref{thm: Main theorem of Neumann diffusions}}
We express the left hand side in matrix form and then compute:
\begin{align*}
    &\sum_{j=1}^{n}(P(X(t) = j\mid X(0) = i) - P(X(t) = j \mid X(0)= k))^2\frac{1}{d_j}\\
    &= (R^t \delta_{i} - R^t \delta_{k})^{\top}T^{-1}_{S}(R^t \delta_{i} - R^t \delta_{k})\\
    &= (R^t(\delta_{i} - \delta_{k}))^{\top}T^{-1}_{S}(R^t(\delta_{i} - \delta_{k}))\\
    &= (\delta_{i} - \delta_{k})^{\top}(R^t)^{\top}T^{-1}_{S}R^t(\delta_{i} - \delta_{k})\\
    &= (\delta_{i} - \delta_{k})^{\top}(C\Sigma^t B^{\top})^{\top}T^{-1}_{S}C\Sigma^t B^{\top}(\delta_{i} - \delta_{k})\\
    &= (\delta_{i} - \delta_{k})^{\top}B\Sigma^t C^{\top}T^{-1}_{S}C\Sigma^t B^{\top}(\delta_{i} - \delta_{k})\\
    &= (\delta_{i} - \delta_{k})^{\top}B\Sigma^t {\underbrace{W^{\top}T^{1/2}T^{-1}T^{1/2}W^{\top}}_I}\Sigma^t\Psi^{\top}(\delta_{i} - \delta_{k})\\
    &= (\Sigma^t B^{\top}(\delta_{i} - \delta_{k})^{\top})^{\top}\Sigma^t B^{\top}(\delta_{i} - \delta_{k})\\
    &= ||g_{t}(i) - g_{t}(k)||^2. 
\end{align*}
\section{Data and code availability}
The data and code for our numerical experiments has been made available at \url{https://github.com/ShashankSule/Neumann_maps/tree/pub}. 

\clearpage
\begin{thebibliography}{00}
\bibitem{SGTbook} Chung, Fan RK. Spectral graph theory. Vol. 92. American Mathematical Soc., 1997.
\bibitem{ChungYau1997} Chung, F.R. and Yau, S.T., 1997. Eigenvalue inequalities for graphs and convex subgraphs. Communications in Analysis and Geometry, 5(4), pp.575-623.
\bibitem{sgt1} Chung, F.R.K., Graham, R.L. and Yau, S.T., 1996. On sampling with Markov chains. Random Structures \& Algorithms, 9(1‐2), pp.55-77.
\bibitem{sgt2} Chung, F.R. and Yau, S.T., 1994. A Harnack inequality for homogeneous graphs and subgraphs. Communications in Analysis and Geometry, 2(4), pp.627-640.
\bibitem{sgt3} Chung, F.R. and Yau, S.T., 1997. Eigenvalue inequalities for graphs and convex subgraphs. Communications in Analysis and Geometry, 5(4), pp.575-623.
\bibitem{sgt4} Tan, J., 1999. Eigenvalue comparison theorems of Neumann laplacian for graphs. Interdisciplinary information sciences, 5(2), pp.157-159.
\bibitem{landmark0} Czaja, W. and Ehler, M., 2012. Schroedinger eigenmaps for the analysis of biomedical data. IEEE Transactions on Pattern Analysis and Machine Intelligence, 35(5), pp.1274-1280.
\bibitem{coifman2005geometric} Coifman, R.R. and Lafon, S., 2006. Diffusion maps. Applied and computational harmonic analysis, 21(1), pp.5-30.
\bibitem{landmark1} De Silva, V. and Tenenbaum, J.B., 2004. Sparse multidimensional scaling using landmark points (Vol. 120). technical report, Stanford University.
\bibitem{landmark2} Silva, V. and Tenenbaum, J., 2002. Global versus local methods in nonlinear dimensionality reduction. Advances in neural information processing systems, 15.
\bibitem{landmark3} Webster, M.A.R.K. and Sheets, H.D., 2010. A practical introduction to landmark-based geometric morphometrics. The paleontological society papers, 16, pp.163-188.
\bibitem{landmark4} Vladymyrov, M. and Carreira-Perpinán, M.Á., 2013. Locally linear landmarks for large-scale manifold learning. In Machine Learning and Knowledge Discovery in Databases: European Conference, ECML PKDD 2013, Prague, Czech Republic, September 23-27, 2013, Proceedings, Part III 13 (pp. 256-271). Springer Berlin Heidelberg.
\bibitem{landmark5} Thongprayoon, C. and Masuda, N., 2024. Online landmark replacement for out-of-sample dimensionality reduction methods. Proceedings of the Royal Society A, 480(2300), p.20230966.
\bibitem{landmark6} Long, A.W. and Ferguson, A.L., 2019. Landmark diffusion maps (L-dMaps): Accelerated manifold learning out-of-sample extension. Applied and Computational Harmonic Analysis, 47(1), pp.190-211.
\bibitem{landmark8} Pai, G., Talmon, R., Bronstein, A. and Kimmel, R., 2019, January. Dimal: Deep isometric manifold learning using sparse geodesic sampling. In 2019 IEEE Winter Conference on Applications of Computer Vision (WACV) (pp. 819-828). IEEE.
\bibitem{kernelreweighting1} Evans, L., Cameron, M.K. and Tiwary, P., 2023. Computing committors in collective variables via Mahalanobis diffusion maps. Applied and Computational Harmonic Analysis, 64, pp.62-101.
\bibitem{kernelreweighting2} Trstanova, Z., Leimkuhler, B. and Lelièvre, T., 2020. Local and global perspectives on diffusion maps in the analysis of molecular systems. Proceedings of the Royal Society A, 476(2233), p.20190036.
\bibitem{kernelreweighting3} Sule, S., Evans, L. and Cameron, M., 2023. Sharp error estimates for target measure diffusion maps with applications to the committor problem. arXiv preprint arXiv:2312.14418.
\bibitem{kernelreweighting4} Hoyos P, Kileel J. Diffusion maps for group-invariant manifolds. arXiv preprint arXiv:2303.16169. 2023 Mar 28.

\bibitem{landmarkdmap1} Yeh, S.Y., Wu, H.T., Talmon, R. and Tsui, M.P., 2024. Landmark Alternating Diffusion. arXiv preprint arXiv:2404.19649.
\bibitem{landmarkdmap2} Shen, C. and Wu, H.T., 2022. Scalability and robustness of spectral embedding: landmark diffusion is all you need. Information and Inference: A Journal of the IMA, 11(4), pp.1527-1595.
\bibitem{landmarkdmap3} Shen, C., Lin, Y.T. and Wu, H.T., 2022. Robust and scalable manifold learning via landmark diffusion for long-term medical signal processing. Journal of Machine Learning Research, 23(86), pp.1-30.
\bibitem{landmarkdmap4} Cahill, N.D., Czaja, W. and Messinger, D.W., 2014, June. Schroedinger eigenmaps with nondiagonal potentials for spatial-spectral clustering of hyperspectral imagery. In Algorithms and technologies for multispectral, hyperspectral, and ultraspectral imagery XX (Vol. 9088, pp. 27-39). SPIE.
\bibitem{landmarkdmap5} Czaja, W., Dong, D., Jabin, P.E. and Ndjakou Njeunje, F.O., 2021. Transport model for feature extraction. SIAM Journal on Mathematics of Data Science, 3(1), pp.321-341.
\bibitem{datafold} Lehmberg et al., (2020). datafold: data-driven models for point clouds and time series on manifolds. Journal of Open Source Software, 5(51), 2283, https://doi.org/10.21105/joss.02283
\bibitem{cvsimportant} Bonati, L., Trizio, E., Rizzi, A. and Parrinello, M., 2023. A unified framework for machine learning collective variables for enhanced sampling simulations: mlcolvar. The Journal of Chemical Physics, 159(1).
\bibitem{cvsdynamical1} Legoll, F. and Lelievre, T., 2010. Effective dynamics using conditional expectations. Nonlinearity, 23(9), p.2131.
\bibitem{cvsdynamical2} Legoll, F., Lelievre, T. and Olla, S., 2017. Pathwise estimates for an effective dynamics. Stochastic Processes and their Applications, 127(9), pp.2841-2863.
\bibitem{cvsdynamical3} Legoll, F. and Lelievre, T., 2011, August. Some remarks on free energy and coarse-graining. In Numerical Analysis of Multiscale Computations: Proceedings of a Winter Workshop at the Banff International Research Station 2009 (pp. 279-329). Berlin, Heidelberg: Springer Berlin Heidelberg.
\bibitem{cvsrr1} Palacio-Rodriguez, K. and Pietrucci, F., 2022. Free energy landscapes, diffusion coefficients, and kinetic rates from transition paths. Journal of chemical theory and computation, 18(8), pp.4639-4648.
\bibitem{cvsrr2} Mouaffac, L., Palacio-Rodriguez, K. and Pietrucci, F., 2023. Optimal reaction coordinates and kinetic rates from the projected dynamics of transition paths. Journal of Chemical Theory and Computation, 19(17), pp.5701-5711.
\bibitem{cvssampling} Yang, Y.I., Shao, Q., Zhang, J., Yang, L. and Gao, Y.Q., 2019. Enhanced sampling in molecular dynamics. The Journal of chemical physics, 151(7).

\bibitem{transitionstate} Gao, Y.Q. and Yang, L., 2006. On the enhanced sampling over energy barriers in molecular dynamics simulations. The Journal of chemical physics, 125(11).
\bibitem{slowfast} Legoll, F., Lelièvre, T., Myerscough, K. and Samaey, G., 2020. Parareal computation of stochastic differential equations with time-scale separation: a numerical convergence study. Computing and Visualization in Science, 23, pp.1-18.
\bibitem{anisotropic} Singer, A., Erban, R., Kevrekidis, I.G. and Coifman, R.R., 2009. Detecting intrinsic slow variables in stochastic dynamical systems by anisotropic diffusion maps. Proceedings of the National Academy of Sciences, 106(38), pp.16090-16095.
\bibitem{equivariance} Batzner, S., Musaelian, A., Sun, L., Geiger, M., Mailoa, J.P., Kornbluth, M., Molinari, N., Smidt, T.E. and Kozinsky, B., 2022. E (3)-equivariant graph neural networks for data-efficient and accurate interatomic potentials. Nature communications, 13(1), p.2453.
\bibitem{cv1} Ribeiro, J.M.L., Bravo, P., Wang, Y. and Tiwary, P., 2018. Reweighted autoencoded variational Bayes for enhanced sampling (RAVE). The Journal of chemical physics, 149(7).
\bibitem{cv2} Vani, B.P., Aranganathan, A., Wang, D. and Tiwary, P., 2023. Alphafold2-rave: From sequence to boltzmann ranking. Journal of chemical theory and computation, 19(14), pp.4351-4354.
\bibitem{cv3} Wang, D., Wang, Y., Evans, L. and Tiwary, P., 2024. From latent dynamics to meaningful representations. Journal of Chemical Theory and Computation, 20(9), pp.3503-3513.
\bibitem{cv4} Belkacemi, Z., Gkeka, P., Lelièvre, T. and Stoltz, G., 2021. Chasing collective variables using autoencoders and biased trajectories. Journal of chemical theory and computation, 18(1), pp.59-78.
\bibitem{cv5} Lelièvre, T., Pigeon, T., Stoltz, G. and Zhang, W., 2024. Analyzing multimodal probability measures with autoencoders. The Journal of Physical Chemistry B, 128(11), pp.2607-2631.

\bibitem{cv5} Rogal, J., Schneider, E. and Tuckerman, M.E., 2019. Neural-network-based path collective variables for enhanced sampling of phase transformations. Physical Review Letters, 123(24), p.245701.

\bibitem{fp1} Nadler, B., Lafon, S., Coifman, R.R. and Kevrekidis, I.G., 2006. Diffusion maps, spectral clustering and reaction coordinates of dynamical systems. Applied and Computational Harmonic Analysis, 21(1), pp.113-127.
\bibitem{fp2} Coifman, R.R., Kevrekidis, I.G., Lafon, S., Maggioni, M. and Nadler, B., 2008. Diffusion maps, reduction coordinates, and low dimensional representation of stochastic systems. Multiscale Modeling \& Simulation, 7(2), pp.842-864.

\bibitem{nystrom} Czaja, W., Doster, T. and Halevy, A. "An overview of numerical acceleration techniques for nonlinear dimension reduction." Recent Applications of Harmonic Analysis to Function Spaces, Differential Equations, and Data Science: Novel Methods in Harmonic Analysis, Volume 2 (2017): 797-829.

\end{thebibliography}
\end{document}